\documentclass[12pt, fleqn]{article}
\usepackage{latexsym,amsfonts,amssymb}

\usepackage{amsfonts}
\usepackage{amssymb}
\usepackage{amsbsy}
\usepackage{amsmath}
\usepackage{epsf}
\usepackage{graphicx}

\newtheorem{theo}{Theorem}% [] %chapter]
\newcommand{\bt}{\begin{theo}}
\newcommand{\et}{\end{theo}}
\newcommand{\bd}{\begin{displaymath}}
\newcommand{\ed}{\end{displaymath}}

\newcommand{\be} {\begin{equation}}
\newcommand{\ee} {\end{equation}}
\newcommand{\ba} {\begin{array}}
\newcommand{\ea} {\end{array}}

\newcommand{\p} {\partial}

\begin{document}

\begin{center}
%%\medskip\\
 {\Large \bf New mathematical model\\
for fluid-glucose-albumin  transport \\ in peritoneal dialysis}\\
%%\medskip\\
{\bf Roman Cherniha~$^1 ,^2$  and  Jacek Waniewski~$^\dag$ }

{\it $^1$~Institute of Mathematics,  NAS of Ukraine,\\
Tereshchenkivs'ka Street 3, 01601 Kyiv, Ukraine\\
$^2$~Department  of  Mathematics,
     National University
     `Kyiv-Mohyla Academy', \\ 2 Skovoroda Street,
     Kyiv  04070 ,  Ukraine}\\
\texttt{cherniha@imath.kiev.ua} \\
{\it  $^\dag$~Institute of Biocybernetics and Biomedical
Engineering,
PAS,\\Ks. Trojdena 4, 02 796 Warszawa, Poland}\\
\texttt{jacekwan@ibib.waw.pl}

\end{center}

\begin{abstract}

A mathematical model for fluid transport  in peritoneal dialysis is
constructed. The model is based on a three-component nonlinear
system of two-dimensional partial differential equations for fluid, glucose and albumin transport with the
relevant boundary and initial conditions. Non-constant
steady-state solutions of the model are studied. The
restrictions on the parameters arising in the model are
established with the aim
 to obtain exact formulae for
the non-constant steady-state solutions.   As the result, the exact
formulae for the  fluid fluxes from blood to tissue and  across the tissue were constructed together with two linear
autonomous ODEs for glucose and albumin concentrations. The analytical results were checked for their
applicability for the description of fluid-glucose-albumin  transport  during
peritoneal dialysis.

\end{abstract}



\medskip
\centerline{{\bf 1 Introduction}}

 Peritoneal dialysis is a life
saving treatment for chronic patients with end stage renal disease
%\cite{gok-nolph}
(Gokal R and   Nolph 1994). Dialysis fluid is infused into the peritoneal
cavity, and,  during its dwell there, small metabolites (urea,
creatinine) and other uremic toxins diffuse from blood to the fluid,
and after some time (usually a few hours) are removed together with
the drained fluid. The treatment is repeated continuously. The
peritoneal transport occurs between dialysis fluid in the peritoneal
 cavity and blood passing down capillaries in tissue surrounding the
peritoneal cavity. The capillaries are distributed within the tissue
at different distance from the tissue surface that is in contact with
dialysis fluid. The solutes, which are transported between blood and
dialysis fluid, have to cross two transport barriers: the capillary
wall and a tissue layer. Typically, many solutes are transported from
blood to dialysate, but some solutes that are present in high
concentration in dialysis fluid are transported to blood. This kind
of transport system happens also in other medical treatments, as
local delivery of anticancer medications, and some experimental or
natural physiological phenomena. Mathematical description of these
systems was obtained using partial differential equations based on
the simplification that capillaries are homogeneously distributed
within the tissue
%% \cite{Fles84,wan-et-al99,Wan2002}
(Flessner et al. 1984; Waniewski et al. 1999; Waniewski 2002). Experimental evidence confirmed the good
applicability of such models
%\cite{Fles85}
(Flessner et al. 1985).

Another objective of peritoneal dialysis is to remove excess water
from the patient
% \cite{gok-nolph}.
(Gokal R and   Nolph 1994). This is gained by inducing high
osmotic pressure in dialysis fluid by adding a solute in high
concentration. The most often used solute is glucose. This medical
application of high osmotic pressure is rather unique for peritoneal
dialysis.
Mathematical description of fluid and solute transport between blood and dialysis
fluid in the peritoneal cavity has not been formulated fully yet, in spite of the well
 known basic physical laws for such transport.
  A previous attempt did not result in a satisfactory description,
  and was disproved later on
  % \cite{seames90,Fles94}.
  (Seames et al. 1990;Flessner  1994).
  Recent mathematical, theoretical and numerical studies
   introduced new concepts on peritoneal transport and yielded
   better results for the transport of fluid and osmotic agent
   % \cite{Ch-wa-2005,joanna-et-al-05,joanna-et-al-06,ch-et-al-2007,wan-et-al-2007,wan-et-al-2009}.
   (Cherniha and  Waniewski 2005; Stachowska-Pietka et al. 2005; Stachowska-Pietka et al. 2006;    Cherniha et al. 2007;Waniewski et al.2007, Waniewski et al.2009).
     However, the problem of
     a combined description of osmotic ultrafiltration to the
     peritoneal cavity, absorption of osmotic agent from the
     peritoneal cavity and leak of macromolecules
     (proteins, e.g., albumin) from blood to the peritoneal cavity
     has not been addressed yet, see for example in
     %\cite{Fles,joanna-et-al-07}.
     (Flessner  2001; Stachowska-Pietka et al. 2007). Therefore, we present here
     a new extended model for these phenomena and investigate
     its mathematical structure.

%%\cite{seams}, \cite{Fles}.
%%The present study is aimed on investigation of some basic questions concerning the role of various transport
%%components, as osmotic gradient and hydrostatic pressure gradient, in the general formulation of the model.

The paper is organized as follows. In section 2, a mathematical
model of glucose and albumin transport in peritoneal dialysis is
constructed. In section 3, non-constant  steady-state solutions of
the model are constructed and their properties are investigated.
Moreover, these solutions are tested for the real parameters that
represent clinical treatments of peritoneal dialysis. The results
are compared with those derived by numerical simulations for a
simplified model
%\cite{ch-et-al-2007}.
(Cherniha et al. 2007).
Finally, we present some
conclusions and discussion in the last section.

 %%the hydrostatic pressure and the glucose concentration   in peritoneal dialysis.

%%\section{Mathematical Model}
\medskip
\centerline{\textbf{2. Mathematical model}}

 The mathematical description of
transport processes within the tissue consists in local balance of
fluid volume and solute mass. For incompressible fluid, the change
of volume may occur due to elasticity of the tissue. The fractional
void volume, i.e. the volume occupied by the fluid in the
interstitium (the rest of the tissue being cells and macromolecules
forming interstitium) expressed per one unit volume of the whole
tissue is denoted $\nu(t,x)$, and its time evolution is described by
the following equation:
 \be \label{1.1} {{\partial \nu} \over
{\partial t}} = - {{\partial j_U } \over {\partial x}} + q_U - q_l
\ee where $j_U(t,x)$ is the volumetric fluid flux across the tissue
(ultrafiltration), $q_U(t,x)$ is the density of volumetric fluid
flux from blood to the tissue, while the density of volumetric fluid
flux from  the tissue to the lymphatic vessels $q_l$ (hereafter we
assume that it is a known positive constant, nevertheless it can be
also a function of hydrostatic pressure) produces absorbtion of
solutes from the tissue.

The independent variables $t$ is time, and $x$ is the distance from
the tissue surface in contact with dialysis fluid (flat geometry of
the tissue is here assumed). The solutes, glucose and albumin, are
distributed only within the interstitial fluid, and their
concentrations in this fluid are denoted by $C_G(t,x)$ and
$C_A(t,x)$, respectively. The equation that describes the local
changes of glucose  amount, $\nu C_G$, is as follows:
 \be
\label{1.2} {{\partial (\nu C_G )} \over {\partial t}} =  -
{{\partial j_G } \over {\partial x}} + q_G, \ee
 where $j_G(t,x)$ is
glucose flux through the tissue, and $q_G(t,x)$ is the density of
glucose flux from blood.

Similarly, equation that describes the local changes of albumin
amount, $\nu C_A$, is as follows:
 \be
\label{1.2-a} {{\partial (\alpha \nu C_A )} \over {\partial t}} =  -
{{\partial j_A } \over {\partial x}} + q_A, \ee
 where $j_A(t,x)$ is
albumin  flux through the tissue, and $q_A(t,x)$ is the density of
albumin flux from blood, while the coefficient $\alpha<1$ takes into
account  that a part of the  fractional void volume $\nu$ is only
accessible for albumin because of big size of molecules
%\cite{Fles,joanna-et-al-07}.
  (Flessner  2001; Stachowska-Pietka et al. 2007).
 The flows of fluid and solutes are
described according to linear non-equilibrium thermodynamics.
Osmotic pressure of glucose and oncotic pressure of  albumin  are
described by van't Hoff law, i.e. it is proportional to the relevant
concentrations.

The
% volumetric
fluid  flux across the tissue is generated by hydrostatic,
osmotic and oncotic pressure gradients:
 \be \label{1.3} j_U = - \nu K{{\partial P} \over
{\partial x}} + \sigma _{TG} \nu KRT{{\partial C_G } \over {\partial
x}}+ \sigma _{TA}\gamma \nu KRT{{\partial C_A } \over {\partial x}},
\ee whereas for the density of fluid flux from blood to tissue we
assume that it is generated by the hydrostatic, osmotic and oncotic
pressure differences between blood and tissue:
 \be \label{1.4} q_U = L_p a(P_{B} - P)- L_p a\sigma_{G} RT(C_{GB} -
 C_G)- \gamma L_p a\sigma_{A} RT(C_{AB} - C_A), \ee
where $P(t,x)$ is hydrostatic  pressure.

  The glucose solute flux across the tissue
is composed of diffusive component (proportional to glucose
concentration gradient) and convective component (proportional to
glucose concentration and fluid flux):
\be \label{1.5} j_G  = -
\nu D_G {{\partial C_G } \over {\partial x}} + S_{TG} C_G j_U. \ee
The density of glucose flux from blood to tissue has diffusive
component (proportional to the difference of glucose concentration
in blood, $C_{GB}$, and glucose concentration in tissue, $C_G$),
convective component (proportional to the density of fluid flow
from blood to tissue, $q_U$) and component presenting lymphatic
absorbtion: \be \label{1.6} q_G = p_G a(C_{GB}-C_G )+S_{G}q_U \left(
{(1-F_G)C_{GB}+F_G C_G}\right) - q_lC_G. \ee

In quite similar way, we can construct   equations for the albumin
solute flux across the tissue,$j_A(t,x)$, and the density of albumin
flux from blood to tissue $q_A(t,x)$:
\be \label{1.5-a} j_A  = -
\alpha \nu D_A {{\partial C_A } \over {\partial x}} + S_{TA} C_A
j_U, \ee \be
 \label{1.6-a} q_A = p_A a(C_{AB}-C_A )+S_{A}q_U \left(
{(1-F_A)C_{AB}+F_A C_A}\right) - q_lC_A. \ee

 The coefficients in the above equations are: $K$ --
hydraulic permeability of tissue, $\sigma_{TG}$ and $\sigma_{TA}$ --
the Staverman reflection
 coefficients for
glucose and albumin in tissue, $S_{TG} = 1-\sigma_{TG}$ and $S_{TA}
= 1-\sigma_{TA}$-- sieving coefficients of glucose and albumin in
tissue, $\sigma_{G}$ and $\sigma_{A}$ -- the Staverman reflection
 coefficients for
glucose and albumin in the capillary wall, $S_{G} = 1-\sigma_{G}$
and $S_{A} = 1-\sigma_{A}$-- sieving coefficients of glucose and
albumin in the capillary wall, $R$ -- gas constant, $T$ --
temperature, $L_p$ -- hydraulic permeability of the capillary wall,
$a$ -- density of capillary surface area,  $D_G$ and $D_A$ --
diffusivities  of glucose and albumin  in tissue,  $p_G$ and $p_A$
-- diffusive permeabilities of the capillary wall for glucose and
albumin, $P_B$ - hydrostatic pressure in blood,
   $F_G\leq 1$ and $F_A\leq 1$  -- weighing factors,
   and $\gamma\leq 1$ is a coefficient that recalculates osmotic pressure of albumin
   to the total oncotic pressure exerted by all proteins
   % \cite{joanna-et-al-07}.
     (Stachowska-Pietka et al. 2007).

Equations (\ref{1.1})-(\ref{1.2-a}) together with equations
(\ref{1.3})-(\ref{1.6-a}) for flows form a system of three nonlinear
partial differential equations with three variables: $\nu, P, C_A$,
and $C_G$. Therefore, an additional, constitutive, equation is
necessary, and this is the equation describing how fractional fluid
volume, $\nu$, depends on interstitial pressure, $P$. This
dependence can be established  using data from experimental studies
%\cite{joanna-et-al-06}.
  ( Stachowska-Pietka et al. 2006). It turns out that \be \label{1.7-a} \nu =
F(P), \ee where $F$ is a monotonically non-decreasing bounded
function with the limits: $F \to \nu_{\min}$ if $P \to P_{\min}$ and
$F \to \nu_{\max}$ if $P \to P_{\max}$ (particularly, one may take
$P_{\min}=-\infty, \ P_{\max}=\infty$).
%\be \label{1.7}
% \nu(P) = \nu_{\min}+ \frac
%{\nu_{\max}-\nu_{\min}}{1+\left( \frac{\nu_{\max }- \nu_{\min }}
%{\nu_0 -\nu_{\min}} -1 \right) \exp ( - b(P-P_0)) }, \ee
Here $\nu_{\min}<1 $ and $ \nu_{max}<1 $ are empirically measured
constants. The reader may find possible analytical forms for the
function $F$ in
% \cite{joanna-et-al-06,ch-et-al-2007}.
  (Stachowska-Pietka et al. 2006; Cherniha et al. 2007).

Boundary conditions for a tissue layer of width $L$ impermeable at
$x = L$ and in contact with dialysis fluid at $x = 0$ are as
follows:
 \be \label{1.8}  x = 0: \ P
= P_D, \quad  C_G = C_{GD}, \quad  C_A = C_{AD} \ee
 \be \label{1.8*}  x = L: \
{{\partial P} \over {\partial x}} = 0,
%%%\quad (or\quad P = P_0),
 \quad {{\partial C_G} \over {\partial x}} = 0, \quad {{\partial C_A} \over {\partial x}} = 0.
\ee Initial conditions describe equilibrium within the tissue
without any contact with dialysis fluid at $x = 0$: \be \label{1.9}
t = 0: \ P = P_0, \ C_G = C_{GB}, \ C_A = C_{AB}. \ee

It is easily seen that equations (\ref{1.1})-(\ref{1.7-a}) can be
united into three nonlinear partial differential equations (PDEs)
for finding the hydrostatic pressure $P(t,x)$, the glucose
concentration $C_G$ and the albumin concentration $C_A$. Thus, these
PDEs together with boundary and initial conditions
(\ref{1.8})-(\ref{1.9}) form the nonlinear  boundary-value problem.
Possible values of the  parameters arising in this problem
  can be established from experimental data published in a wide range of papers.
 % (see, e.g., \cite{Fles, Wan2001, Wan2002,joanna-et-al-06,joanna-et-al-07, ch-et-al-2007,Wan2008}.

Finally, we note that finding  the function $j_U(t,x)$, presenting
the fluid fluid flux across the tissue and
calculating ultrafiltration from the tissue to peritoneal cavity, is
the most important.

\newpage

\centerline{\textbf{3. Steady-state solutions of the  model and
their applications}}
%% and their interpretation.}

First of all we consider the special case: if the boundary
conditions (\ref{1.8}) is replaced by zero Neumann conditions at
$x=0$ then the steady state solution can be easily found because it
does't depend on $x$. In fact solving algebraic equations \be
\label{2.0} q_U - q_l =0, \quad q_G =0, \quad q_A =0, \ee one easily
obtains the constant  steady-state solution \be \label{2.*}
\ba{l}\medskip C_G^*= \frac{p_G a+q_lS_{G}(1-F_G)}{p_G
a+q_l(1-S_{G}F_G)} C_{GB}\\ \medskip C_A^*= \frac{p_A
a+q_lS_{A}(1-F_A)}{p_A a+q_l(1-S_{A}F_A)}
C_{AB}\\
P^* = P_B -q_l\Bigl(\frac{1}{L_pa}+RT \Bigl(\frac{\sigma_G^2}{p_G
a+q_l(1-S_{G}F_G)}+\frac{\gamma\sigma_A^2}{p_A
a+q_l(1-S_{A}F_A)}\Bigr)\Bigr).\ea \ee In the case $q_l=0$, i.e.,
zero flux from the tissue to the lymphatic vessels, formulae
(\ref{2.*}) produce \be \label{2.**} C_G^*= C_{GB}, \quad  C_A^*=
C_{AB}, \quad P^* = P_B,  \ee otherwise \be \label{2.***} C_G^*<
C_{GB}, \quad  C_A^* < C_{AB}, \quad P^* < P_B.  \ee
 However, such
solution cannot describe any peritoneal transport.

To find non-constant  steady-state solutions,  we reduce Eqs.
(\ref{1.1})-(\ref{1.2-a}) to an equivalent form by introducing
non-dimensional independent and dependent variables \be \label{2.1}
%%\ba{l}
x^*=\frac xL, \quad t^*=\frac{K(P_D-P_0)t}{L^2}, \ee \be \label{2.2}
p(t^*,x^*)= \frac{P-P_0}{P_D-P_0}, \quad u(t^*,x^*) =
\frac{C_G-C_{GB}}{C_{GD}-G_{GB}}, \quad
w(t^*,x^*)=\frac{C_A}{C_{GD}-G_{GB}}.
 \ee
Thus, after rather simple calculations and taking into account
 Eq. (\ref{1.3}),  (\ref{1.5}), and (\ref{1.5-a}), one obtains
 Eq. (\ref{1.1})-(\ref{1.2-a}) in  the forms (hereafter
 upper index $*$
is omitted) \be \label{2.3} \frac{1}{t_0}\frac{\p\nu}{\p t}=
\frac{1}{t_0}\frac{\p}{\p x} \Bigl({\nu{\partial p} \over {\partial
x}}\Bigr) - \sigma _{1} \frac{\p}{\p x} \Bigl({\nu{\partial u} \over
{\partial x}}\Bigr)- \sigma _{2} \frac{\p}{\p x} \Bigl({\nu{\partial
w} \over {\partial x}}\Bigr)  + q_U - q_l,\ee
\medskip
 \[ \frac{1}{t_0}\frac{\p(\nu u)}{\p t} + \frac{\sigma_{TG}u_0 }{t_0}\frac{\p\nu}{\p t}=
d _{1} \frac{\p}{\p x} \Bigl({\nu{\partial u} \over {\partial
x}}\Bigr)+\frac{S_{TG}}{t_0}\frac{\p}{\p x} \Bigl({u\nu{\partial p}
\over {\partial x}}\Bigr) - S_{TG}\sigma _{1} \frac{\p}{\p x}
\Bigl({u\nu{\partial u} \over {\partial x}}\Bigr)\] \be \label{2.4}
- S_{TG}\sigma _{2} \frac{\p}{\p x} \Bigl({u\nu{\partial w} \over
{\partial x}}\Bigr)+(u_0(S_G-S_{TG})+S_GF_Gu)q_U
-b_1u-\sigma_{TG}u_0q_l,\ee\medskip
 \[ \frac{1}{t_0}\frac{\p(\nu w^*)}{\p t} + \frac{ \sigma_{TA}w_0}{t_0} \frac{\p\nu}{\p t}=
d _{2} \frac{\p}{\p x} \Bigl({\nu{\partial w} \over {\partial
x}}\Bigr)+\frac{S_{TA}}{t_0}\frac{\p}{\p x} \Bigl({w^*\nu{\partial
p} \over {\partial x}}\Bigr) - S_{TA}\sigma _{1} \frac{\p}{\p x}
\Bigl({w^*\nu{\partial u} \over {\partial x}}\Bigr)\] \be
\label{2.5} - S_{TA}\sigma _{2} \frac{\p}{\p x}
\Bigl({w^*\nu{\partial w} \over {\partial
x}}\Bigr)+(w_0(S_A-S_{TA})+S_AF_Aw^*)q_U
-b_2w^*-\sigma_{TA}w_0q_l,\ee
 where \be \label{2.6} \ba{l}
 \medskip
q_U= \frac{L_paL^2}{K}\Bigl(
\frac{1}{t_0}(p_0-p)+\frac{\sigma_{G}\sigma_{1}}{\sigma_{TG}}u+
\frac{\sigma_{A}\sigma_{2}}{\sigma_{TA}}w^*
\Bigr),\\
  \sigma_1=\sigma_{TG}KRT \frac{C_{GD}-G_{GB}}{L^2}, \quad
  \sigma_2=\sigma_{TA}KRT \gamma \frac{C_{GD}-G_{GB}}{L^2},\\
d_1=\frac{D_{G}}{L^2}, \quad d_2=\frac{\alpha D_{A}}{L^2}, \\
%s_1=\frac{S_{TG}}{L}, \quad s_2=\frac{S_{TA}}{L},\\
 b_1= p_Ga+q_l, \quad  b_2= p_Aa+q_l, \\
 u_0=\frac{C_{GB}}{C_{GD}-G_{GB}},   \quad
 w_0=\frac{C_{AB}}{C_{GD}-G_{GB}}, \quad p_0 = \frac{P_B-P_0}{P_D-P_0}\\

 %%h_1=L_pa RT\sigma_{CG}, b_0=P_GaC_{GB}-S_{CG}(1-F_G)h_0C_{GB},\\
t_0 = \frac{L^2}{K(P_D-P_0)},\quad w^*=w-w_0, \quad
%%b_1=P_Ga-(1-2F_G)S_Gh_0, \quad b_2=S_{CG}F_Gh_1.
 \ea \ee

We want to find  the steady state solutions of Eqs.
(\ref{2.3})-(\ref{2.5}) satisfying the boundary conditions
(\ref{1.8})-(\ref{1.8*}). They take the form \be \label{2.7} x = 0:
\ p = 1, \quad u = 1, \quad  w = \frac{C_{AD}}{C_{GD}-G_{GB}} \ee
 \be \label{2.8}  x = 1: \
{{\partial p} \over {\partial x}} = 0,
 \quad {{\partial u} \over {\partial x}} = 0, \quad {{\partial w} \over {\partial x}} = 0.
\ee for the non-dimensional variables.

\medskip

Obviously, Eqs. (\ref{2.3})-(\ref{2.5}) can be  reduced to the
ordinary differential equations (ODEs) to find steady state
solutions: \be \label{2.10}  \frac{1}{t_0}\frac{d}{d x} \Bigl({\nu{d
p} \over {d x}}\Bigr) - \sigma _{1} \frac{d}{d x} \Bigl({\nu{d u}
\over {d x}}\Bigr)- \sigma _{2} \frac{d}{d x} \Bigl({\nu{d w} \over
{d x}}\Bigr) + q_U - q_l = 0,\ee
\medskip
 \[
d _{1} \frac{d}{d x} \Bigl({\nu{d u} \over {d
x}}\Bigr)+\frac{S_{TG}}{t_0}\frac{d}{d x} \Bigl({u\nu{d p} \over {d
x}}\Bigr) - S_{TG}\sigma _{1} \frac{d}{d x} \Bigl({u\nu{d u} \over
{d x}}\Bigr)\] \be \label{2.11} - S_{TG}\sigma _{2} \frac{d}{d x}
\Bigl({u\nu{d w} \over {d x}}\Bigr)+(u_0(S_G-S_{TG})+S_GF_Gu)q_U
-b_1u-\sigma_{TG}u_0q_l =0 ,\ee\medskip
 \[
d _{2} \frac{d}{d x} \Bigl({\nu{d w} \over {d
x}}\Bigr)+\frac{S_{TA}}{t_0}\frac{d}{d x} \Bigl({w^*\nu{d p} \over
{d x}}\Bigr) - S_{TA}\sigma _{1} \frac{d}{d x} \Bigl({w^*\nu{d u}
\over {d x}}\Bigr)\] \be \label{2.12} - S_{TA}\sigma _{2} \frac{d}{d
x} \Bigl({w^*\nu{d w} \over {d
x}}\Bigr)+(w_0(S_A-S_{TA})+S_AF_Aw^*)q_U -b_2w^*-\sigma_{TA}w_0q_l =
0.\ee  Non-linear system of ODEs (\ref{2.10})-(\ref{2.12}) is still
very complicated for integrations with arbitrary coefficients. Thus,
we look for the correctly-specified coefficients when this system
can be simplified. It can be noted that the relations \be
\label{2.13} S_A=S_{TA}, \quad S_G=S_{TG} \ee lead to an essential
simplification of this system. In fact, using equation (\ref{2.10}),
expressions for $q_U$ from (\ref{2.6}) and   $j_U$ from (\ref{1.3}),
rewritten in non-dimensional variables
 \be \label{2.13a} j_U =  L\nu \Bigl(-\frac{1}{t_0}{{\partial p} \over
{\partial x}} + \sigma _{1}  {{\partial u } \over {\partial x}}+
\sigma _{2} {{\partial w } \over {\partial x}}\Bigr), \ee
 one arrives
at the semicoupled system of ODEs \be \label{2.14}
\frac{K}{L_paL^2}\frac{d}{d x}\Bigl({\nu{d q_U} \over {d x}}\Bigr)= q_U - q_l \ee
\be \label{2.15} j_U = \frac{K\nu}{L_paL}\frac{dq_U}{d x} \ee to
find functions $q_U$ and $j_U$.

Equation  (\ref{2.14}) is the linear second-order ODE provided the
function $\nu$ is known. Nevertheless $\nu$ depends on the pressure,
which is  also to be found function, we may use additional
restrictions on the given function $F$ from (\ref{1.7-a}). For
example, assuming that $F$ is proportional to the inverse function
to $P(x)$, we obtain \be \label{2.16} \nu(x) = \nu_m, \ee where
$\nu_m$ is a positive constant. Substituting (\ref{2.16}) into
system  (\ref{2.14})--(\ref{2.15}),  we easily  find its general
solution: \be \label{2.17} q_U = C_1e^{-\lambda x} + C_2e^{\lambda
x} + q_l, \ee \be \label{2.18} j_U = \frac L\lambda (-C_1e^{-\lambda
x} + C_2e^{\lambda x}),  \quad  \lambda=
\sqrt{\frac{L_paL^2}{K\nu_m}}.\ee The arbitrary constants $C_1$ and
$C_2$  can be specified using the boundary conditions
(\ref{2.7})--(\ref{2.8}) since the functions $q_U$ and $j_U$ are
expressed via $p, \ u, \ w$ and their first-order derivatives (see
formulae (\ref{2.6}) and (\ref{2.13a})). Making rather simple
calculations, one obtains \be \label{2.19} C_1 = (q_0
-q_l)\frac{e^{2\lambda}}{1+ e^{2\lambda}}, \quad C_2 = (q_0
-q_l)\frac{1}{1+ e^{2\lambda}},\ee where \be \label{2.20} q_0 =
\frac{L_paL^2}{K}\Bigl(
\frac{1}{t_0}(p_0-1)+\sigma_1+\sigma_{2}\frac{C_{AD}-C_{AB}}{C_{GD}-G_{GB}}\Bigr).\ee

Having the explicit  formulae for $q_U$ and $j_U$, system of ODEs
(\ref{2.11}) and (\ref{2.12}) with restrictions (\ref{2.13}) can be
reduced to two linear autonomous ODEs \be \label{2.21} d _{1}\nu_m
\frac{d^2u}{d x^2} +\frac{S_{G}}{\lambda}(C_1e^{-\lambda x} -
C_2e^{\lambda x})\frac{du}{d x} +\Bigl( f_1(C_1e^{-\lambda x} +
C_2e^{\lambda x})-\kappa_1\Bigr)u-u_{01} =0 \ee and \be \label{2.22}
d _{2}\nu_m \frac{d^2w}{d x^2} +\frac{S_{A}}{\lambda}(C_1e^{-\lambda
x} - C_2e^{\lambda x})\frac{dw}{d x} +\Bigl( f_2(C_1e^{-\lambda x} +
C_2e^{\lambda x})-\kappa_2\Bigr)(w-w_0)-w_{01} =0 \ee
 to find the functions
$u(x)$ and $w(x)$. Hereafter the notations  \be \label{2.23} \ba{l}
f_1=S_GF_G-S_G, \quad \kappa_1=p_Ga+(1-S_GF_G)q_l, \quad
u_{01}=\sigma_{G}u_0q_l, \\  f_2=S_AF_A-S_A, \quad
\kappa_2=p_Aa+(1-S_AF_A)q_l, \quad w_{01}=\sigma_{A}w_0q_l \ea \ee
are used. To our best knowledge, the general solutions of ODEs
(\ref{2.21}) and (\ref{2.22}) (obviously, both equations have the
same structure) are unknown in explicit form. On the other hand,
since the unknown functions $u(x)$ and $w(x)$ should satisfy the
boundary conditions (\ref{2.7})--(\ref{2.8}), the corresponding
linear problems can be numerically solved using, for example, Maple
program package. Finally, using two expressions for $q_U$ from
(\ref{2.6}) and (\ref{2.17}), we obtain   the function   \be
\label{2.24}p(x)=
 p_0+t_0\sigma_{1}u+t_0\sigma_{2}(w-w_0)-\frac{t_0K}{L_paL^2}\Bigl( C_1e^{-\lambda x} + C_2e^{\lambda x} -
q_l\Bigr). \ee

In the next section, the realistic  values of parameters arising in
the formulae derived above  will be established  and used to
calculate the steady-state solutions of the  model with restrictions
(\ref{2.13}) and (\ref{2.16}).

\vspace{5mm}

Since restriction (\ref{2.16}) is rather an artificial
simplification for modeling   fluid transport in peritoneal dialysis
and it  needs additional justifications,  we examined the cases when
the function $\nu$ is non-constant and satisfies the properties
presented after formula (\ref{1.7-a}). The simplest case occurs when
$\nu$ is linear monotonically decreasing function \be \label{2.25}
\nu(p(x))\equiv \nu(x)=\nu_{max} -(\nu_{max}-\nu_{min})x, \quad x
\in [0,1]. \ee Substituting (\ref{2.25}) into (\ref{2.13}), we
obtain the linear ODE \be \label{2.26} (\nu_{max}
-(\nu_{max}-\nu_{min})x){{d^2 q_U} \over {d x^2}}  -
(\nu_{max}-\nu_{min}){{d q_U} \over {d x}}- \frac{L_paL^2}{K}(q_U -
q_l) =0.\ee It turns out that this ODE reduces to the modified
Bessel equation of the zero order
(see, e.g.,
%\cite{bateman}
 Bateman 1974) \be \label{2.27}
y^2\frac{d^2q_V}{dy^2} + y\frac{dq_V}{dy} -y^2q_V = 0 \ee by the
substitution \be \label{2.28} y^2 = 4\delta_*(\nu_* - x), \ q_V =
q_U - q_l, \ \nu_*=\frac{\nu_{max}}{\nu_{max}-\nu_{min}}>1, \
%\quad
\delta_*= \frac{L_paL^2}{K(\nu_{max}-\nu_{min})}>0\ee The general
solution of (\ref{2.27}) is well-known, hence, using formulae
(\ref{2.28}), one obtains the solution of (\ref{2.26}): \be
\label{2.29} q_U = C_1 I_0(2\sqrt{\delta_*(\nu_*-x)}) +
C_2K_0(2\sqrt{\delta_*(\nu_*-x)}) + q_l,\ee where $I_0$ and $K_0$
are the modified Bessel functions of the first and third kind,
respectively.

Substituting the function $q_U$ obtained into (\ref{2.14}) and using
the well-known relations between the Bessel functions
%\cite{bateman}
( Bateman 1974), we find the function: \be \label{2.30} j_U
=-L\sqrt{\frac{\nu_*-x}{\delta_*}}\Bigl( C_1
I_1(2\sqrt{\delta_*(\nu_*-x)}) -
C_2K_1(2\sqrt{\delta_*(\nu_*-x)})\Bigr),\ee where $I_1$ and $K_1$
are the modified Bessel functions of the first order. The arbitrary
constants $C_1$ and $C_2$  can be specified  in the quite similar
way as this was done in the case $\nu(x)=const$. Omitting  rather
simple calculations, we present only the result: \be \label{2.31}
C_1 = \frac{(q_0-q_l)K_1(2\sqrt{\delta_*(\nu_*-1)})}
{I_0(2\sqrt{\delta_*\nu_*})K_1(2\sqrt{\delta_*(\nu_*-1)})+
K_0(2\sqrt{\delta_*\nu_*})I_1(2\sqrt{\delta_*(\nu_*-1)})}, \ee \be
\label{2.32} C_2 = \frac{(q_0-q_l)I_1(2\sqrt{\delta_*(\nu_*-1)})}
{I_0(2\sqrt{\delta_*\nu_*})K_1(2\sqrt{\delta_*(\nu_*-1)})+
K_0(2\sqrt{\delta_*\nu_*})I_1(2\sqrt{\delta_*(\nu_*-1)})}, \ee where
$q_0$ is defined by (\ref{2.20}).

Thus, we have found  the explicit  formulae  for  $q_U$ and $j_U$.
In contrast to the case  $\nu(x)=const$,  they contain not only
elementary functions but  transcendental functions as well. Having
formulae (\ref{2.29})--(\ref{2.30}),   system of ODEs (\ref{2.11})
and (\ref{2.12}) with restrictions (\ref{2.13})  can be reduced to
two linear autonomous ODEs to find the functions $u(x)$ and $w(x)$.
These equations possess the forms
%\[
\be \label{2.33}\ba{l} d_{1}(\nu_{max}-\nu_{min})
\Bigl((\nu_*-x)\frac{d^2u}{d x^2}-\frac{du}{d x} \Bigr) \\
-\frac{S_{G}}{\sqrt\delta_*}\frac{d}{d x}\Bigl(\sqrt{\nu_*-x}(C_1
I_1(2\sqrt{\delta_*(\nu_*-x)}) -
C_2K_1(2\sqrt{\delta_*(\nu_*-x)}))u\Bigr)\\ +\Bigl( S_GF_G(C_1
I_0(2\sqrt{\delta_*(\nu_*-x)}) +
C_2K_0(2\sqrt{\delta_*(\nu_*-x)}))-\kappa_1\Bigr)u-u_{01} =0 \ea\ee

and

\be \label{2.34} \ba{l} d_{2}(\nu_{max}-\nu_{min})
\Bigl((\nu_*-x)\frac{d^2w}{d x^2}-\frac{dw}{d x} \Bigr) \\
-\frac{S_{A}}{\sqrt\delta_*}\frac{d}{d x}\Bigl(\sqrt{\nu_*-x}(C_1
I_1(2\sqrt{\delta_*(\nu_*-x)}) -
C_2K_1(2\sqrt{\delta_*(\nu_*-x)}))(w-w_0)\Bigr) \ea \ee \[+\Bigl(
S_AF_A(C_1 I_0(2\sqrt{\delta_*(\nu_*-x)}) +
C_2K_0(2\sqrt{\delta_*(\nu_*-x)}))-\kappa_2\Bigr)(w-w_0)-w_{01} =0\]

Nevertheless both equations are  linear second order ODEs with the
same structure, it seems to be unrealistic to construct their
general solutions because of their awkwardness. Thus, we will
numerically solve them together with the boundary conditions
(\ref{2.7})--(\ref{2.8}), using Maple program package. In the next
section, the realistic values of the parameters arising in formulae
(\ref{2.29})--(\ref{2.34}) will be established  and used to
calculate the steady-state solutions.

\centerline{\textbf{4. Applications for the peritoneal transport}}

Here we present the results of application of the formulae derived
in Section 3. Our aim is
to check
 whether they are
applicable for describing the fluid-glucose-albumin  transport  in
peritoneal dialysis. The parameters arising in the formulae were derived from experimental data and
applied in previous mathematical studies
% \cite{ch-et-al-2007,Fles, Wan2001, Zaka,Wan2004,Imho}.
(Imholz et al. 1998;  Zakaria et al. 1999;  Flessner  2001; Waniewski 2001; Waniewski 2004; Stachowska-Pietka et al. 2006; Cherniha et al. 2007).

Thus, we used the following values of parameters and absolute
constants:\\
$K=5.14 \cdot 10^{-5}$ $cm^{2} \cdot min^{-1} \cdot mmHg^{-1}$ -
hydraulic permeability of tissue, $RT=18 \cdot 10^{3}$ $mmHg \cdot
mmol^{-1} \cdot mL$ - gas constant times temperature, $L=1.0$ $cm$ -
width of the  tissue layer, $L_{P}a=7.3 \cdot 10^{-5}$ $mL\cdot
min^{-1} \cdot mmHg^{-1} \cdot g^{-1}$- hydraulic permeability of
the capillary wall, $q_l =0.26\cdot 10^{-4}$ $mL \cdot min^{-1}
\cdot cm^{-3}$ - volumetric fluid flux from  the tissue to the
lymphatic vessels,
 $D_G=12.11 \cdot 10^{-5}$ $cm^{2} \cdot min^{-1}$ -
diffusivity of glucose in tissue divided by $\nu_{0}$, $D_A=0.2
\cdot 10^{-5}$ $cm^{2} \cdot min^{-1}$ - diffusivity of albumin in
tissue divided by $\nu_{0}$, $p_{G}a=3.4 \cdot 10^{-2}$ $mL \cdot
min^{-1} \cdot g^{-1}$ - diffusive permeability of the capillary
wall for glucose, $p_{A}a=
3\cdot 10^{-4}$ $mL \cdot min^{-1} \cdot g^{-1}$ - diffusive
permeability of the capillary wall for albumin (there are no precise values
for this parameter in literature, however, $p_{A}$ is expected to be
at least 100 times smaller than $p_{G}$ \cite{joanna-et-al-07} ),
$\sigma_{TG}$ - the Staverman reflection coefficient for glucose in
tissue (varies from $0$ to $0.01$), $S_{TG} = 1-\sigma_{TG}$ -
sieving coefficient of glucose in tissue, $\sigma_{TA}$ - the
Staverman reflection coefficient for albumin in tissue (varies from
$0.05$ to $0.5$), $S_{TA} = 1-\sigma_{TA}$ - sieving coefficient of
albumin in tissue,
 $C_{GB}=6\cdot 10^{-3}$ $mmol \cdot mL^{-1}$ and
$C_{AB}=0.6\cdot 10^{-3}$ $mmol \cdot mL^{-1}$-glucose and albumin
concentration in the blood, respectively,  $C_{GD}=180\cdot 10^{-3}$
$mmol\cdot mL^{-1}$ and $C_{AD}=0$ - glucose and albumin
concentration in the dialysate, respectively, $P_{B}=15$ $mmHg$  -
hydrostatic pressure in the blood, $P_{D}=12$ $mmHg$ -
intraperitoneal hydrostatic pressure, $P_{0}=0$  - initial
interstitial hydrostatic pressure, $F_G=0.5$ and $F_A=0.5$  --
weighing factors for glucose and
 albumin, respectively, and the non-dimensional coefficients
$\nu_{\min}=0.17, \nu_{max}=0.35, \nu_0=0.17$, $\alpha=0.8, \gamma
=1$.

Let us consider the first case, when  restrictions (\ref{2.13}) and
(\ref{2.16}) take place. First of all, it seems to be reasonable to
set $\nu_m= (\nu_{max}+\nu_{min})/2 = 0.26,$ i.e., we assume that
the fractional void volume at the  steady state stage of the the
peritoneal transport is an intermediate value between its maximum and
minimum. The Staverman reflection coefficients for glucose
and albumin in tissue
 are now $\sigma_{TG}=0.001$ and $\sigma_{TA}=0.25$,
respectively.

\begin{figure}[t]
\begin{minipage}[t]{5.5cm}
\centerline{\includegraphics[width=6.0cm]{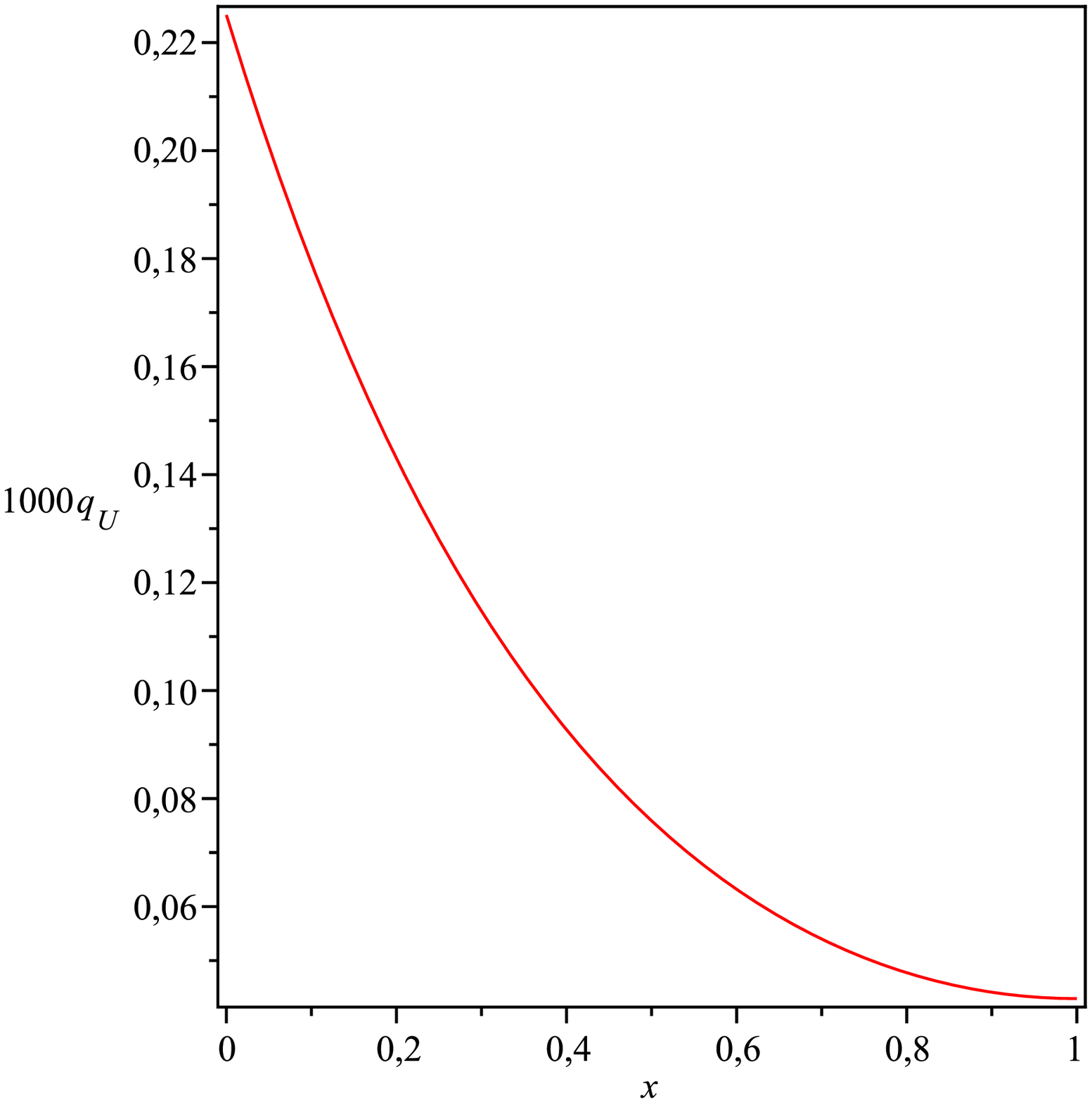}}
\end{minipage}
\hfill
\begin{minipage}[t]{5.5cm}
\centerline{\includegraphics[width=6.0cm]{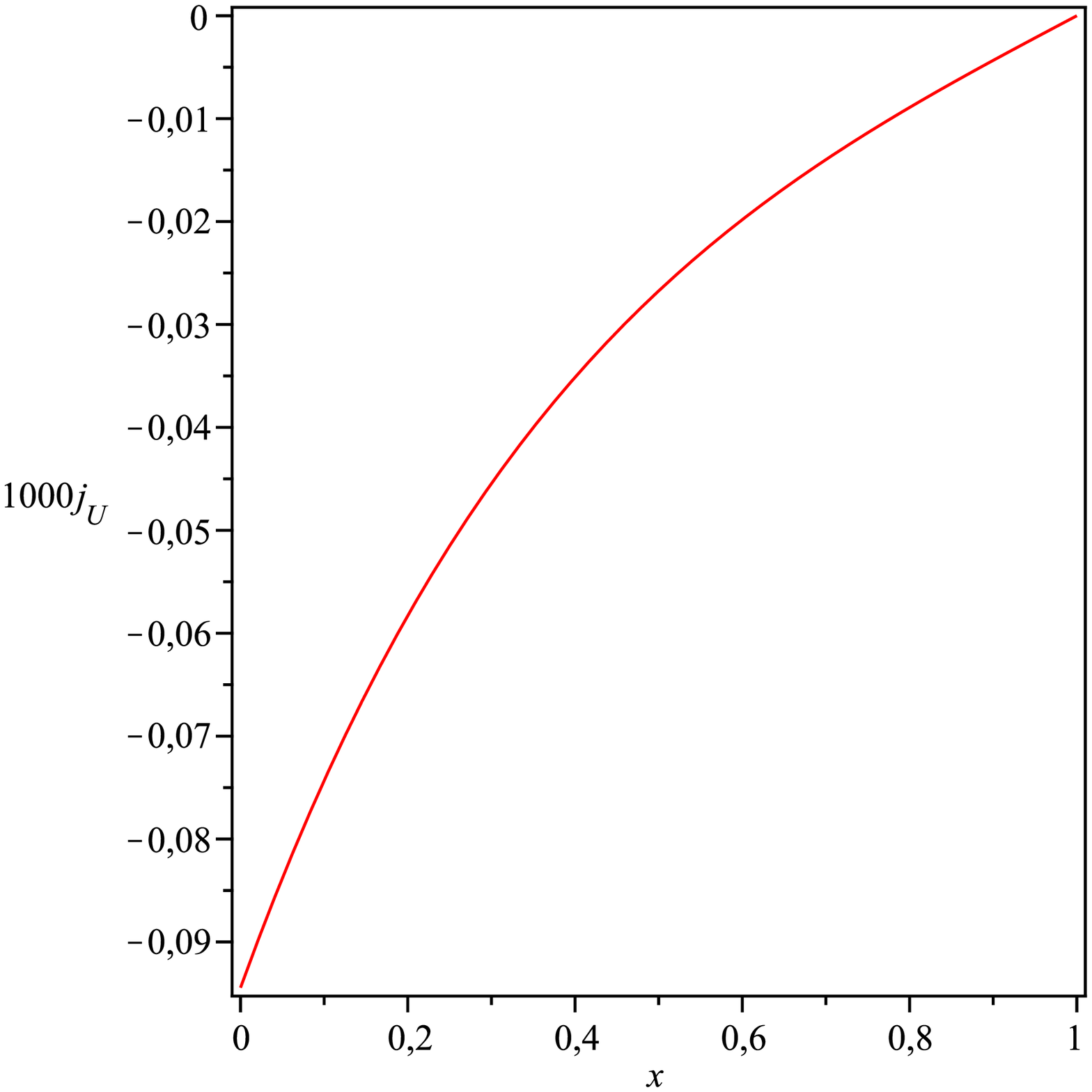}}
\end{minipage}
\caption{The  fluid flux from blood to tissue $q_{U}$ (in $mL \cdot min^{-1}
\cdot cm^{-3}$) and  the  fluid flux  across  tissue
 $j_{U}$ (in $mL \cdot min^{-1}
\cdot cm^{-2}$ )  in the case $\nu= (\nu_{max}+\nu_{min})/2 = 0.26$,
$\sigma_{TG}=0.001$,  and $\sigma_{TA}=0.25$.}
\label{Fig-1}
 \end{figure}

Fig.1 presents the space distributions of the density of fluid flux
from blood to tissue $q_U$ and the fluid flux across  tissue $j_U$, calculated
using formulae (\ref{2.17})--(\ref{2.20}).    One sees that the
function $q_U(x)$ is monotonically decreasing, while the function
$j_U(x)$ is monotonically increasing with the distance from the peritoneal surface, and this corresponds to the
experimental data and numerical simulations for the simplified model
%\cite{ch-et-al-2007,wan-et-al-2007}.
(Cherniha et al. 2007; Waniewski et al. 2007).

\begin{figure}[t]
\begin{minipage}[t]{5.5cm}
\centerline{\includegraphics[width=6.0cm]{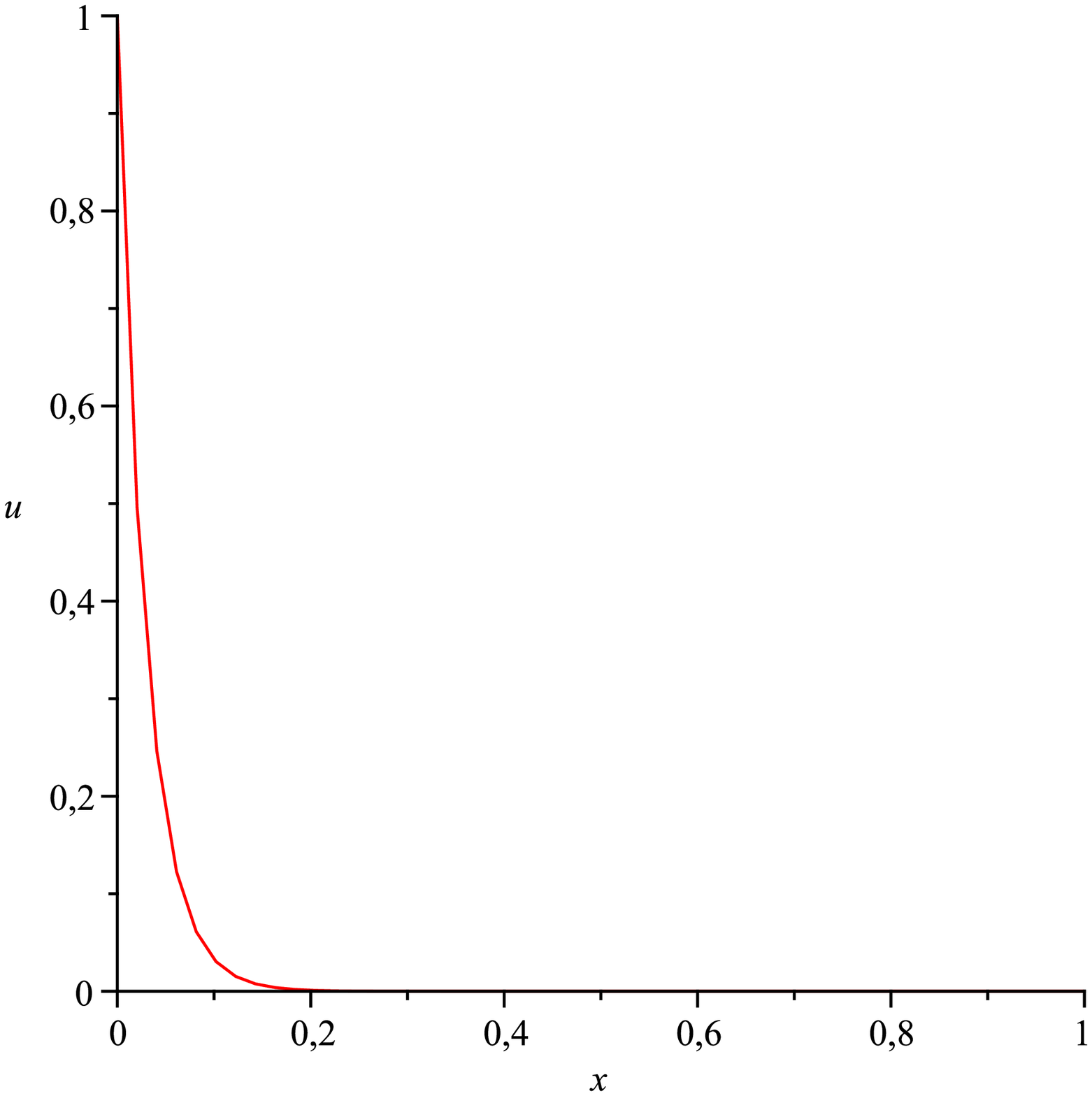}}
\end{minipage}
\hfill
\begin{minipage}[t]{5.5cm}
\centerline{\includegraphics[width=6.0cm]{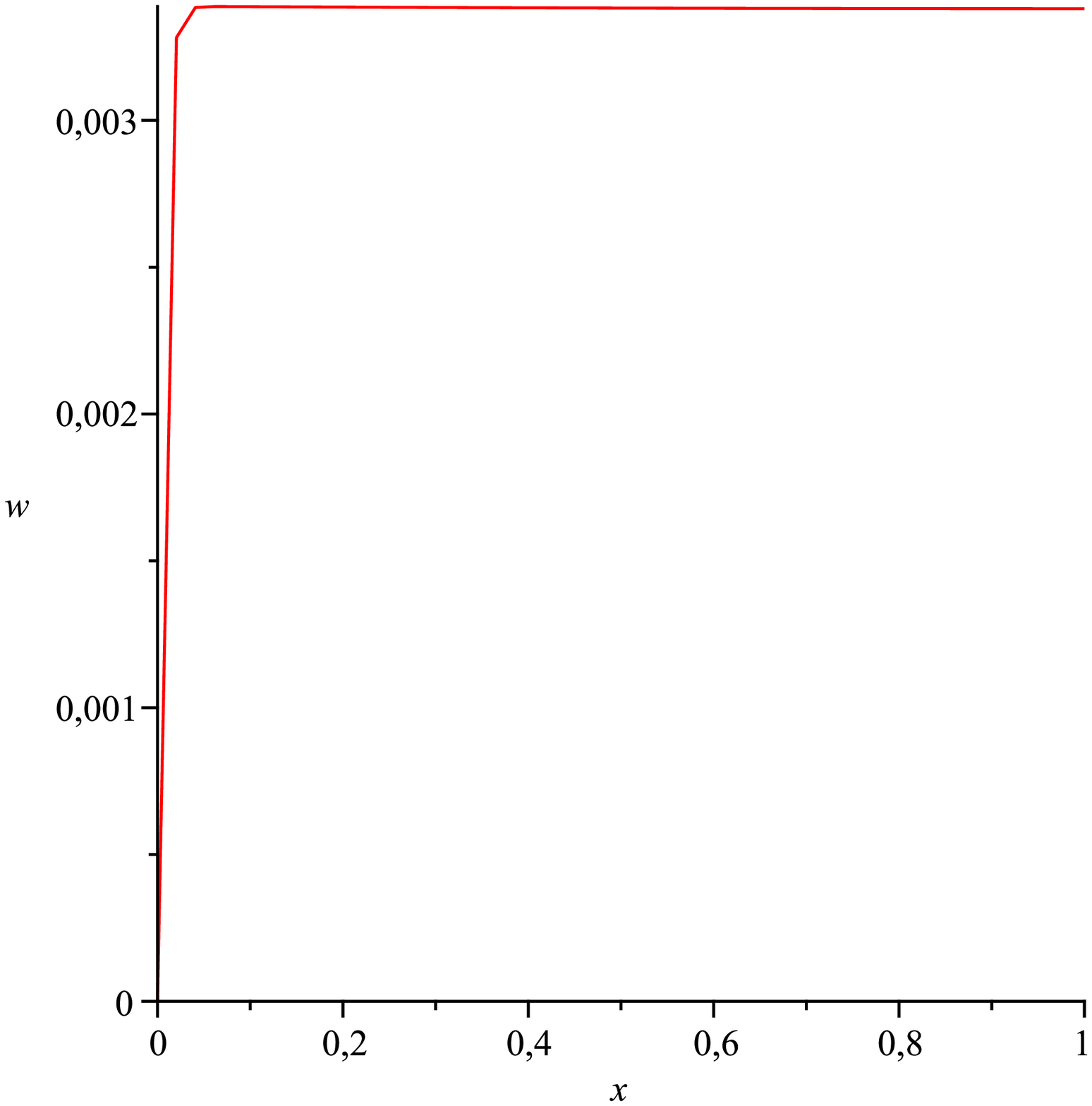}}
\end{minipage}
\caption{The  non-dimensional glucose concentration $u$ and albumin concentration $w$  in the case $\nu= (\nu_{max}+\nu_{min})/2 = 0.26$,
$\sigma_{TG}=0.001$,  $\sigma_{TA}=0.25$, and  $p_{A}a=
3\cdot 10^{-4}$ $mL \cdot min^{-1} \cdot g^{-1}.$  }
\label{Fig-2}
 \end{figure}

% \newpage

 \begin{figure}[t]
\begin{minipage}[t]{5.5cm}
\centerline{\includegraphics[width=6.0cm]{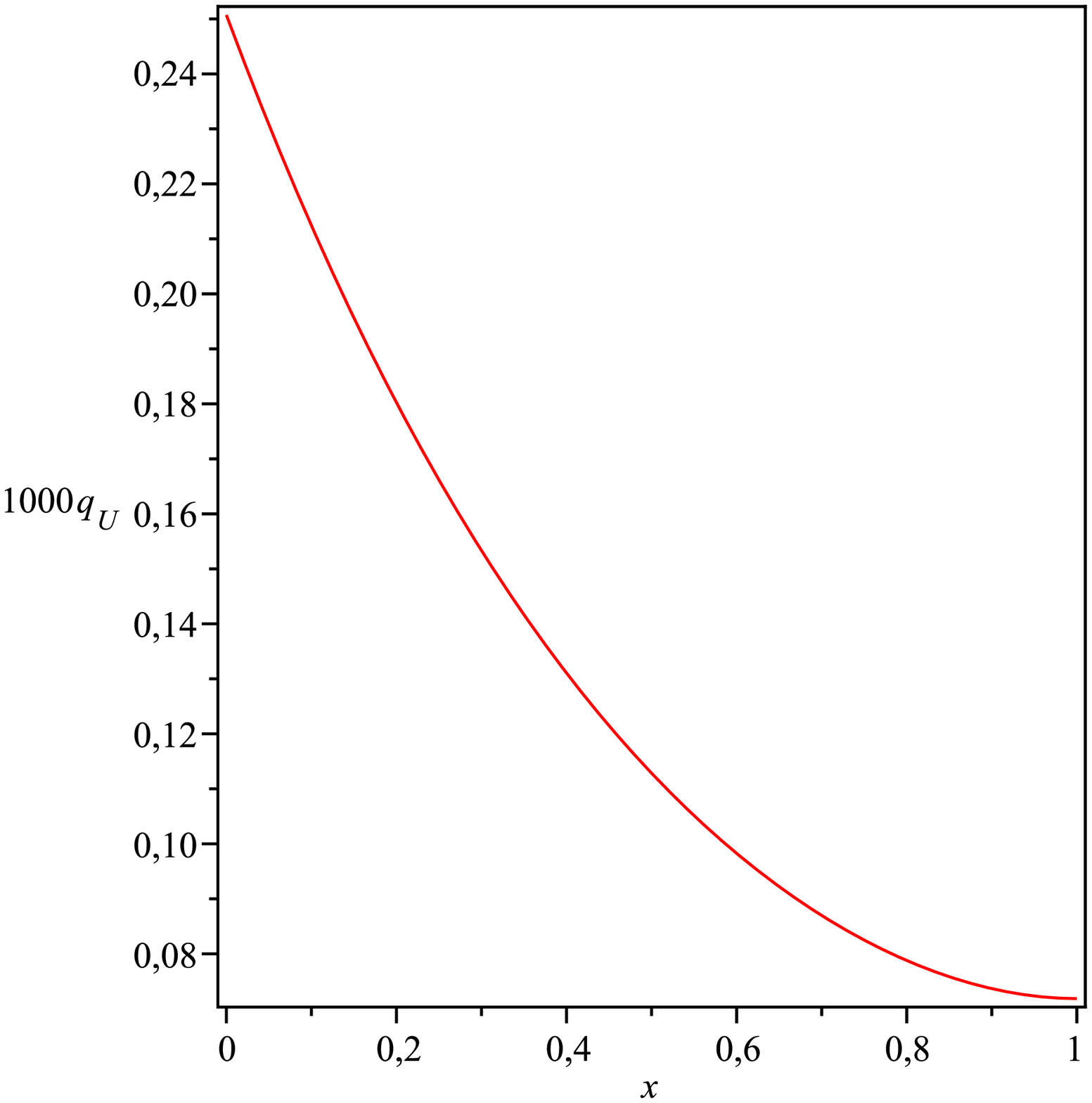}}
\end{minipage}
\hfill
\begin{minipage}[t]{5.5cm}
\centerline{\includegraphics[width=6.0cm]{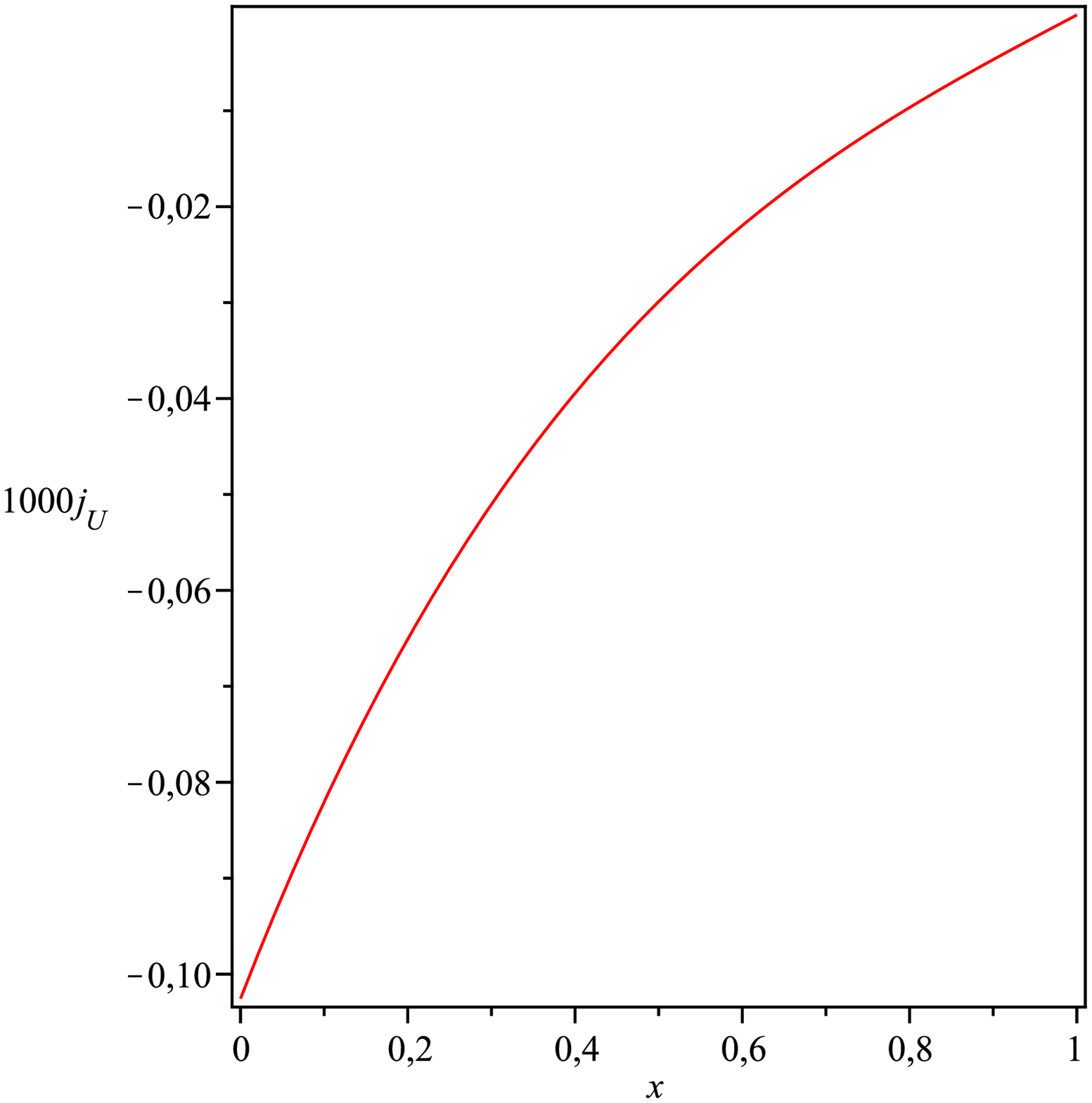}}
\end{minipage}
\caption{The  fluid flux from blood to tissue $q_{U}$ (in $mL \cdot min^{-1}
\cdot cm^{-3}$) and  the  fluid flux  across  tissue
 $j_{U}$ (in $mL \cdot min^{-1}
\cdot cm^{-2}$ )  in the case $\nu=\nu_{max} -(\nu_{max}-\nu_{min})x$,
$\sigma_{TG}=0.001$,  and $\sigma_{TA}=0.25$.}
\label{Fig-3}
 \end{figure}

\begin{figure}[t]
\begin{minipage}[t]{5.5cm}
\centerline{\includegraphics[width=6.0cm]{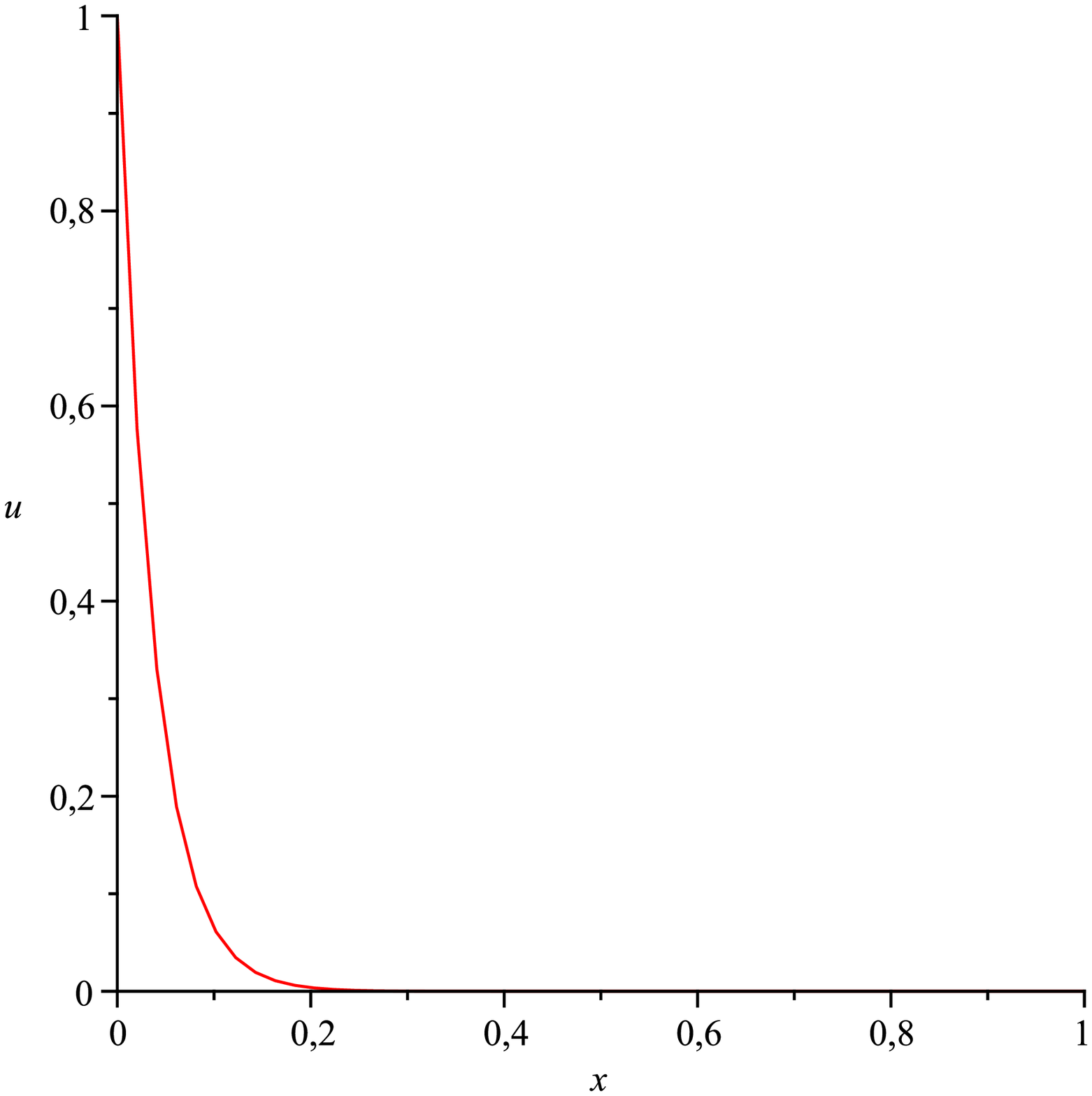}}
\end{minipage}
\hfill
\begin{minipage}[t]{5.5cm}
\centerline{\includegraphics[width=6.0cm]{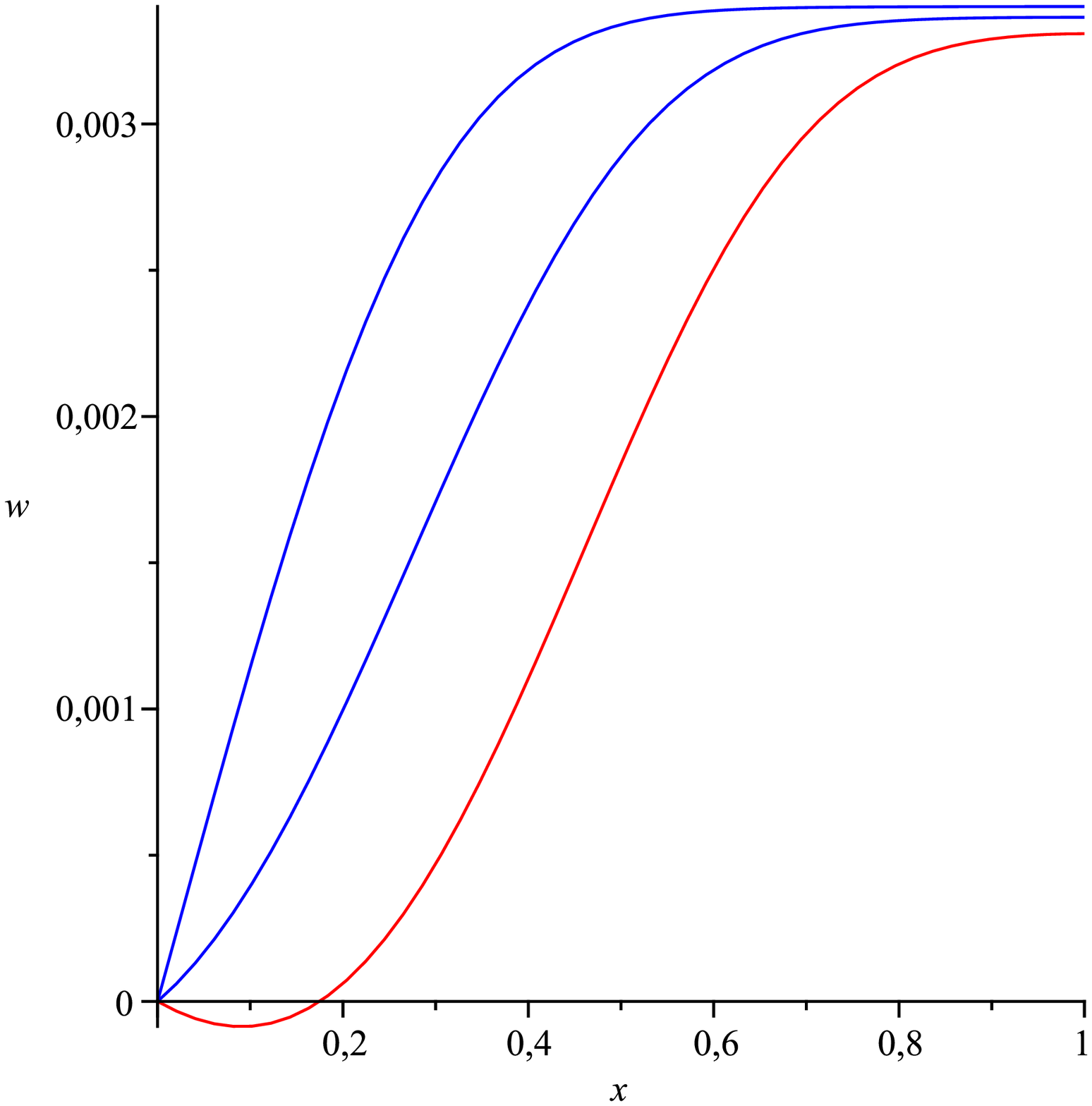}}
\end{minipage}
\caption{The  non-dimensional glucose concentration $u$ and albumin concentration $w$  in the case $\nu=\nu_{max} -(\nu_{max}-\nu_{min})x$,
$\sigma_{TG}=0.001$,  $\sigma_{TA}=0.25$, and $p_{A}a=
5\cdot 10^{-4},$    $ 3\cdot 10^{-4},$   $ 2\cdot 10^{-4}$ $mL \cdot min^{-1} \cdot g^{-1}$  }
\label{Fig-4}
 \end{figure}

Using the value of  the fluid flux $j_U$ at the point $x=0$,
one may calculate the reverse water flow (i.e. out of the tissue to
the cavity). Total fluid outflow from the tissue to the cavity,
calculated assuming that the surface area of the contact between
dialysis fluid and peritoneum is equal to $0.5$ $m^{2}$, is about
$0.50$ $mL/min$. Note the similar value was obtained in
(Cherniha et al. 2007) using numerical simulations for the simplified
model.

Fig.2 presents the space distributions of the glucose and albumin
concentrations. The glucose  concentration $u$ decreases rapidly
 with the distance from the peritoneal cavity to zero in
the deeper tissue layer and is
practically zero for $x>0.15 $. As one may expect, the albumin
concentration $w$ increases  with the distance from the peritoneal
cavity. However, the concentration increases rapidly only in a very
thin layer, while for any $x>0.05 $ it is practically constant,
$w\simeq 0.0034$, and in agreement with the initial albumin profile
$C_A=C_{AB}=0.6\cdot 10^{-3}$ $mmol \cdot mL^{-1}$ (see
(\ref{1.9})). However, such sharp increase in concentration is rather unrealistic
and is a consequence of the assumption (\ref{2.16}).

Let us consider now the second  case, which is  more realistic.
We assume that  restrictions (\ref{2.13}) and (\ref{2.25})
take place. The Staverman reflection coefficients were chosen the
same as in the first case, i.e., $\sigma_{TG}=0.001$ and $\sigma_{TA}=0.25$. The results
are presented on Fig. 3 and 4. It is quite interesting that the
profiles for the functions $q_U(x)$ and $j_U(x)$ pictured on Fig.3
are very similar to those on Fig. 1, although the relevant
formulae  are essentially different (the reader may compare
(\ref{2.29}) -- (\ref{2.32}) with (\ref{2.17})--(\ref{2.20})).
Moreover, the total fluid outflow from the tissue to the cavity,
calculated under the above mentioned assumption    is approximately
equal to $-0.55$ $mL/min$, that is only about 10 percents higher than in the previous case.

The space distributions of the glucose and albumin concentrations
are pictured on Fig.4.  The glucose  concentration $u$ is again a
decreasing function, however the tissue layer with non-vanishing $u$
is wider.

The main difference occurs in the case of the albumin
concentration $w$. We found that the albumin concentration $w$  essentially
depends on the parameter $p_{A}a$.  Three curves pictured on Fig.4
correspond to the values  of the  diffusive
permeability  $p_{A}a=
5\cdot 10^{-4},$    $ 3\cdot 10^{-4},$  and  $ 2\cdot 10^{-4}$ $mL \cdot min^{-1} \cdot g^{-1}$,
respectively. Negative values
occurring close to the peritoneal cavity (the  red curve corresponding  to the smallest $p_{A}a$
value) are very small and can be  either a consequence of the simplifications in the model
or an error in numerical solving ODE (\ref{2.34}). We may  interpret these negative values  as follows: there is no albumin
in this layer of the tissue because it was already removed to the
peritoneal cavity.

 Finally, one observes that the relevant curves on Fig.2 and Fig.4
(the middle  curve) are essentially different, nevertheless they were
obtained for the same parameters. As follows from the Fig.4, the
albumin concentration $w$ increases in deeper layers  of tissue as
well, and is equal to  the constant, corresponding to the initial
profile, only close to the opposite side of tissue. Thus, such
steady-state profile of the  albumin concentration  is more
realistic than in the first case (see Fig.2).

\medskip

%% All pressures are scaled assuming atmospheric pressure equal to zero.

\centerline{\textbf{5. Conclusions}}
%% and their interpretation.}
\medskip
In this paper, a new  mathematical model for fluid transport  in
peritoneal dialysis was constructed. The model is based on a
three-component nonlinear system of two-dimensional partial
differential equations and the relevant boundary and initial
conditions. To analyze the non-constant steady-state solutions, the
model was reduced to the non-dimensional form. Having in mind to
obtain exact formulae for such solutions, we found the restrictions
on parameters arising in the model that essentially simplified the
equations of the model. As the result, the exact formulae for the  density of
fluid flux from blood to tissue  and the fluid flux across the
tissue were constructed together with two linear autonomous  ODEs for
glucose and albumin concentrations.

The analytical results were checked, whether they are applicable for
describing the fluid-glucose-albumin  transport  in peritoneal dialysis.
Thus, the realistic values of the parameters, arising in the
formulae, were used to calculate the steady-state solutions
of the model. The conclusion is rather optimistic
because even the simplest approximation of fractional fluid volume
$\nu$ via linear function with the  correctly specified coefficients
leads to plausible results. Other, more realistic approximation of $\nu$ may in future result in similar exact formulae.
 However, in general the assumption about equality of the reflection coefficients in the tissue and in the capillary wall, although it demonstrates an interesting specific symmetry in the equations, can be too restrictive for practical applications of the derived formulae
  %(c.f. ref. 13).
  ( Waniewski et al. 2009). Therefore, other approaches to find  the analytical solutions of the model need to be looked for.
%steady-state solutions.

%\newpage

 \centerline{\textbf{  Acknowledgments.}}
\medskip
This work was done within the  project "Mathematical models of fluid
and solute transport in normal  and pathological  tissue" supported
by Mianowski Fund (Warsaw, Poland). R.Ch. thanks the Department of
Mathematical Modelling of Physiological Processes (IBIB of PAS,
Warsaw), where the main part of this work was carried out, for
hospitality.

\renewcommand{\refname}{References}

\end{document}